\RequirePackage{ifpdf}
\ifpdf 
\documentclass[pdftex]{sigma}
\else
\documentclass{sigma}
\fi

\begin{document}

\allowdisplaybreaks

\renewcommand{\thefootnote}{$\star$}

\renewcommand{\PaperNumber}{093}

\FirstPageHeading

\newcommand{\mR}{\mathbb{R}}
\newcommand{\mC}{\mathbb{C}}
\newcommand{\mN}{\mathbb{N}}
\newcommand{\mE}{\mathbb{E}}
\newcommand{\mZ}{\mathbb{Z}}
\newcommand{\mS}{\mathbb{S}}
\newcommand{\mP}{\mathbb{P}}

\newcommand{\cD}{\mathcal{D}}
\newcommand{\cM}{\mathcal{M}}
\newcommand{\cH}{\mathcal{H}}
\newcommand{\cE}{\mathcal{E}}
\newcommand{\cF}{\mathcal{F}}
\newcommand{\cP}{\mathcal{P}}
\newcommand{\cC}{\mathcal{C}}
\newcommand{\cR}{\mathcal{R}}
\newcommand{\cK}{\mathcal{K}}
\newcommand{\cS}{\mathcal{S}}

\newcommand{\ux}{\underline{x}}
\newcommand{\uxb}{\underline{x} \grave{}}
\newcommand{\uxi}{\underline{\xi}}
\newcommand{\uyb}{\underline{y} \grave{}}
\newcommand{\uy}{\underline{y}}
\newcommand{\uub}{\underline{u} \grave{}}
\newcommand{\uu}{\underline{u}}
\newcommand{\uvb}{\underline{v} \grave{}}
\newcommand{\uv}{\underline{v}}

\newcommand{\ds}{d \sigma_{\ux}}
\newcommand{\dsb}{d \sigma_{\uxb}}
\newcommand{\dsx}{d \sigma_{x}}

\newcommand{\pj}{\partial_{x_j}}
\newcommand{\pjb}{\partial_{{x \grave{}}_{j}}}
\newcommand{\pkb}{\partial_{{x \grave{}}_{k}}}

\newcommand{\pI}{\partial_{x_i}}
\newcommand{\pIb}{\partial_{{x \grave{}}_{i}}}
\newcommand{\pIcont}{\partial_{x_i} \rfloor}
\newcommand{\pIbcont}{\partial_{{x \grave{}}_{i}} \rfloor}

\newcommand{\pk}{\partial_{x_k}}
\newcommand{\px}{\partial_x}
\newcommand{\py}{\partial_y}
\newcommand{\upx}{\partial_{\underline{x}}}
\newcommand{\upxb}{\partial_{\underline{{x \grave{}}} }}
\newcommand{\upy}{\partial_{\underline{y}}}
\newcommand{\upyb}{\partial_{\underline{{y \grave{}}} }}

\ShortArticleName{Hermite Polynomials Related to the Dunkl Laplacian}

\ArticleName{An Alternative Def\/inition of the Hermite Polynomials\\ Related to the Dunkl Laplacian\footnote{This paper is a contribution to the Special
Issue on Dunkl Operators and Related Topics. The full collection
is available at
\href{http://www.emis.de/journals/SIGMA/Dunkl_operators.html}{http://www.emis.de/journals/SIGMA/Dunkl\_{}operators.html}}}

\Author{Hendrik DE BIE}

\AuthorNameForHeading{H. De Bie}

\Address{Department of Mathematical Analysis, Faculty of Engineering, Ghent University,\\ Krijgslaan 281, 9000 Gent, Belgium}

\Email{\href{mailto:Hendrik.DeBie@UGent.be}{Hendrik.DeBie@UGent.be}}

\URLaddress{\url{http://cage.ugent.be/~hdebie/}}

\ArticleDates{Received October 07, 2008, in f\/inal form December 18,
2008; Published online December 28, 2008}

\Abstract{We introduce the so-called Clif\/ford--Hermite polynomials in the framework of Dunkl operators, based on the theory of Clif\/ford analysis. Several properties of these polynomials are obtained, such as a Rodrigues formula, a dif\/ferential equation and an explicit relation connecting them with the generalized Laguerre polynomials. A link is established with the generalized Hermite polynomials related to the Dunkl operators (see [R\"osler M., \textit{Comm. Math. Phys.} {\bf 192} (1998), 519--542, \href{http://arxiv.org/abs/q-alg/9703006}{q-alg/9703006}]) as well as with the basis of the weighted $L^{2}$ space introduced by Dunkl.}

\Keywords{Hermite polynomials; Dunkl operators; Clif\/ford analysis}

\Classification{33C80; 33C45; 30G35}

\section{Introduction}

Dunkl operators (see \cite{MR951883, MR1827871}) are combinations of dif\/ferential and dif\/ference operators, associated to a f\/inite ref\/lection group $G$. One of the interesting aspects of these operators is that they allow for the construction of a Dunkl Laplacian, which is a combination of the classical Laplacian in $\mR^m$ with some dif\/ference terms, such that the resulting operator is only invariant under $G$ and not under the whole orthogonal group. Moreover, they are directly related to quantum integrable models of Calogero type (see e.g.~\cite{CMS}) and have as such received a lot of attention in the physics literature.

In \cite{MR1620515} generalizations of the classical Hermite polynomials to the framework of Dunkl ope\-ra\-tors were introduced and some of their properties proven, such as a dif\/ferential equation, a~Mehler formula, etc. However, the precise form of these generalized Hermite polynomials is not very clear. Only for special choices of the group~$G$ it is possible to obtain more detailed information and to relate them to classical orthogonal polynomials.

It is the aim of the present paper to introduce generalized Hermite polynomials in a dif\/ferent way. We draw inspiration from Clif\/ford analysis (see a.o.~\cite{MR697564, MR1169463}), a function theory for the Dirac operator, e.g.~in $\mR^m$. Generalizations of the Hermite polynomials to this framework are called Clif\/ford--Hermite polynomials and were introduced by Sommen in \cite{MR926831}. Detailed accounts can also be found in \cite{MR1169463} and in \cite{DBS3} for the superspace case (this can be seen as the study of dif\/ferential operators invariant under the action of the group $O(m)\times Sp(2n)$).
In this paper, we adapt the def\/inition of the Clif\/ford--Hermite polynomials to the case of Dunkl operators.
The advantage of these polynomials is that their relation with classical orthogonal polynomials on the real line is established, leading to a much more concrete form for them. We are also able to prove that the Clif\/ford--Hermite polynomials in the Dunkl framework coincide with the Hermite polynomials introduced by R\"osler in \cite{MR1620515} if a suitable basis of the space of homogeneous polynomials is chosen. We also point out that they have been introduced previously in a dif\/ferent fashion by Dunkl as a basis for a weighted $L^{2}$ space (see \cite{Du4}) and that they have been studied more thoroughly in~\cite{Said}.

Another interesting aspect is that the Clif\/ford--Hermite polynomials always take the same form, be it in the case of the classical Laplace operator with invariance $O(m)$, the case of the Dunkl Laplacian with invariance $G \subset O(m)$ or the super Laplace operator with invariance $O(m)\times Sp(2n)$. The only dif\/ference is a numerical parameter on which these polynomials depend; this parameter can be interpreted as the dimension of the associated theory.

The paper is organized as follows. In Section~\ref{section2} we f\/irst discuss how the Clif\/ford--Hermite polynomials arise in Clif\/ford analysis. Then we give some background on Dunkl operators and we give the def\/inition of the generalized Hermite polynomials due to R\"osler. In Section~\ref{CHdunkl} we introduce the Clif\/ford--Hermite polynomials related to the Dunkl Laplacian. Some basic properties, such as a Rodrigues formula and a dif\/ferential equation are proven. We obtain an expression of the Clif\/ford--Hermite polynomials in terms of the generalized Laguerre polynomials. Then we show that the Clif\/ford--Hermite polynomials and the generalized Hermite polynomials generate the same eigenspace of the associated dif\/ferential-dif\/ference operator and that they coincide if a suitable basis of the space of homogeneous polynomials is chosen.

\section{Preliminaries}\label{section2}

\subsection[Hermite polynomials in Clifford analysis]{Hermite polynomials in Clif\/ford analysis}
\label{hermCA}

The basic operators in harmonic analysis in $\mR^{m}$ are the Laplacian $\Delta = \sum\limits_{i=1}^{m} \partial_{x_{i}}^{2}$, the Euler operator  $\mE = \sum\limits_{i=1}^{m} x_i \partial_{x_i}$ and the norm squared of a vector $x \in\mR^m$, $|x|^2= \sum\limits_{i=1}^{m} x_i^2 $. It is easy to check that the operators
\[ E:= \tfrac{1}{2}|x|^2,\qquad F:= -\tfrac{1}{2}\Delta \qquad\text{and}\qquad
H:= \mE +m/2
\]
satisfy the def\/ining relations of the Lie algebra $\mathfrak{sl}_2$.
These relations are given by
\begin{gather*}
\big[H,E\big] = 2E,\qquad \big[H,F\big] = -2F,\qquad \big[E,F\big] = H.
\end{gather*}

It is now possible to introduce a ref\/inement of harmonic analysis in $\mR^{m}$ by introducing Clif\/ford algebra elements. If we denote by $\cC l_{0,m}$ the orthogonal Clif\/ford algebra of signature $(-1,\ldots, -1)$ generated by $m$ generators $e_{i}$ satisfying
\[
e_{i} e_{j} + e_{j} e_{i} = -2 \delta_{ij},
\]
we can introduce the Dirac operator in $\mR^{m}$ as
\begin{gather}
\upx = \sum_{i=1}^{m} e_{i} \partial_{x_{i}}.
\label{diracoperator}
\end{gather}
We can furthermore identify a vector $x$ in $\mR^{m}$ with the Clif\/ford algebra valued element $\ux$ given~by
\[
\ux = \sum_{i=1}^{m} e_{i} x_{i}.
\]
It is easy to check that $\upx^{2} = - \Delta$ and $\ux^{2} = - |x|^2$.
Moreover, a short calculation shows that
\[
\left\{ \ux, \upx \right\} = \ux \upx + \upx \ux =- ( 2 \mE + m).
\]
Using these relations, it can be checked that the operators $\upx$, $\ux$ together with the previous three operators $E$, $F$ and $H$ generate a Lie superalgebra isomorphic with $\mathfrak{osp}(1|2)$. The even part of this superalgebra being  $\mathfrak{sl}_2$, we have in fact obtained a ref\/inement of harmonic analysis. The function theory related to the Dirac operator (\ref{diracoperator}) is called Clif\/ford analysis (see a.o.~\cite{MR697564, MR1169463, MR1130821} and references therein).

A central notion of this function theory is that of a monogenic function. A $C^{1}$-function $f$, def\/ined in an open set $\Omega \subset \mR^{m}$, with values in the Clif\/ford algebra $\cC l_{0,m}$ is called monogenic if $\upx f = 0$ in $\Omega$. The concept of monogenicity is a higher dimensional analogue of holomorphicity in the complex plane and monogenic functions satisfy several properties similar to those of holomorphic functions, such as Cauchy integral formulae, Taylor and Laurent expansions, \ldots (for a nice overview, see e.g.~\cite{MR1169463}).

As monogenic functions lie at the heart of Clif\/ford analysis, it is important to have tools available to construct such functions. Two of the most important techniques are Fueter's theorem (see a.o.~\cite{MR1509515, MR1812650}) and Cauchy-Kowaleskaia (CK) or monogenic extension (see e.g.~\cite{MR697564, MR0807257}). Fueter's theorem allows to construct monogenic functions starting from a holomorphic function in the complex plane. CK-extension on the other hand allows under certain conditions to construct a monogenic function $F$ in $\mR^{m+1}$, starting from an analytic function $f(\ux)$ in $\mR^{m}$. This extension is given by the formula
\begin{gather}
F(\ux, x_{m+1}) = \sum_{k=0}^{\infty} \frac{x_{m+1}^{k}}{k!} (e_{m+1} \upx)^{k} f(\ux).
\label{CKext}
\end{gather}
It is an easy exercise to check that indeed $(\upx + e_{m+1} \partial_{x_{m+1}}) F = 0$ and that the restriction to~$\mR^{m}$ of $F(\ux, x_{m+1})$ equals $f(\ux)$.

Further, we denote by $\cP = \mR[x_1,\ldots,x_m]$ the space of polynomials on $\mR^m$. We have that $\cP = \oplus_{k=0}^{\infty} \cP_k$ with $\cP_k$ the space of homogeneous polynomials of degree $k$. A spherical monogenic of degree $k$ is an element $M_{k} \in \cP_{k }\otimes \cC l_{0,m}$ satisfying $\upx M_{k} =0$. Spherical monogenics play a~role in Clif\/ford analysis, similar to that of spherical harmonics in harmonic analysis.

It is possible to introduce multi-variable orthogonal polynomials in Clif\/ford analysis. A~gene\-ralization of e.g.\ the Hermite polynomials was f\/irst introduced by Sommen in \cite{MR926831} using the technique of CK or monogenic extension. Recall that the classical Hermite polynomials on the real line are given by the generating function
\[
e^{2 tx - t^{2}} = \sum_{n=0}^{\infty} H_{n}(x) \frac{t^{n}}{n!}
\]
which can be rewritten as
\begin{gather}
e^{z^{2}} = \sum_{n=0}^{\infty} e^{x^{2}} H_{n}(i x) \frac{t^{n}}{n!}, \qquad z = x + it.
\label{gaussext}
\end{gather}
In other words, the Hermite polynomials appear in a natural way when calculating the holomorphic extension of the function $e^{x^{2}}$ in a special way. Similarly, in Clif\/ford analysis generalized Hermite polynomials appear when calculating the monogenic extension of the Gaussian
\[
e^{-|x|^2}
\]
in $\mR^{m}$ to $\mR^{m+1}$. In fact, it is even more general to consider the monogenic extension of
\begin{gather}
e^{-|x|^2} M_{k}
\label{monogenicext}
\end{gather}
with $M_{k}$ a spherical monogenic of a certain degree $k$. This is still in correspondence with formula~(\ref{gaussext}), as the one-dimensional Dirac operator only has the constants as polynomial null-solutions.

Calculating the monogenic extension of (\ref{monogenicext}) using formula (\ref{CKext}) then leads to the following def\/inition (see \cite{MR926831}).

\begin{definition}
Let $M_k$ be a spherical monogenic of degree $k$ and $t$ a positive integer. Then
\[
C\!H_{t}^{m}(M_k) = (D_+)^t M_k,
\]
with
\[
D_{+} = -\upx+2\ux
\]
is a Clif\/ford--Hermite polynomial of degree $t$ associated with $M_{k}$.
\end{definition}

Using this def\/inition we see that the f\/irst few Clif\/ford--Hermite polynomials have the following explicit form:
\begin{gather*}
C\!H_{0}^{m}(M_k) =M_k,\\
C\!H_{1}^{m}(M_k) =2 \ux M_k,\\
C\!H_{2}^{m}(M_k) = \left(4 \ux^2 + 2(2k+m) \right)M_k,\\
C\!H_{3}^{m}(M_k) =\left(8 \ux^3 + 4(2k+m+2)\ux \right)M_k,\\
C\!H_{4}^{m}(M_k) =\left(16 \ux^4 + 16 (2k+m+2)\ux^2+4(2k+m+2)(2k+m)\right) M_k.
\end{gather*}
It is important to note that in order to calculate these examples, only the $\mathfrak{osp}(1|2)$ action is needed and contains all the information. This will allow us to transfer this def\/inition to the Dunkl case. For more properties of the Clif\/ford--Hermite polynomials we refer the reader to~\cite{MR926831, MR1169463, DBS3}.

Furthermore, if we calculate the square of $D_{+}$ we obtain
\begin{gather*}
D_{+}^{2} = (-\upx+2\ux)^{2}
= -\Delta + 4 \ux^{2} - 2\left\{ \upx, \ux \right\}\\
\phantom{D_{+}^{2}}{} =-\Delta - 4 |x|^2 + 4\mE + 2m.
\end{gather*}
As this operator is scalar, it makes sense to let it act on a spherical harmonic instead of on a spherical monogenic. We will use a similar operator in Section~\ref{CHdunkl} to def\/ine Clif\/ford--Hermite polynomials in the setting of Dunkl operators.

It is interesting to note that the Clif\/ford--Hermite polynomials have several applications, e.g.\ in the theory of wavelets (see \cite{MR1755145, NDS} and references therein).

Recently, also an extension of harmonic analysis and Clif\/ford analysis to so-called superspaces has been proposed (see a.o.~\cite{DBS5, DBS3, HDB1}). Superspaces are spaces which are equipped not only with a set of commuting variables, but also with a set of $2n$ anti-commuting variables. It is possible to extend the $\mathfrak{sl}_{2}$-relations to this case. This means deforming $\Delta$, $|x|^{2}$ and $\mE$ to operators which are invariant under $O(m) \times Sp(2n)$ (note that the full symmetry is in fact given by the Lie superalgebra $\mathfrak{osp}(m|2n)$, but that is not relevant for our discussion). The symplectic group has a natural action on the Grassmann algebra generated by the anti-commuting variables. In this framework also, the Clif\/ford--Hermite polynomials appear in a natural way and are used to describe harmonic oscillators and eigenfunctions of a generalized Fourier transform (see~\cite{DBS3, HDB1}).

\subsection{Dunkl operators and Hermite polynomials}

Denote by $\langle \cdot,\cdot \rangle$ the standard Euclidean scalar product in $\mR^{m}$ and by $|x| = \langle x, x\rangle^{1/2}$ the associated norm. For $\alpha \in \mR^{m} - \{ 0\}$, the ref\/lection $r_{\alpha}$ in the hyperplane orthogonal to $\alpha$ is given by
\[
r_{\alpha}(x) = x - 2 \frac{\langle \alpha, x\rangle}{|\alpha|^{2}}\alpha, \qquad x \in \mR^{m}.
\]

A root system is a f\/inite subset $R \subset \mR^{m}$ of non-zero vectors such that, for every $\alpha \in R$, the associated ref\/lection $r_{\alpha}$ preserves $R$. We will assume that $R$ is reduced, i.e.\ $R \cap \mR \alpha = \{ \pm \alpha\}$ for all $\alpha \in R$. Each root system can be written as a disjoint union $R = R_{+} \cup (-R_{+})$, where $R_{+}$ and~$-R_{+}$ are separated by a hyperplane through the origin. The subgroup $G \subset O(m)$ generated by the ref\/lections $\{r_{\alpha} | \alpha \in R\}$ is called the f\/inite ref\/lection group associated with $R$. We will also assume that $R$ is normalized such that $\langle \alpha, \alpha\rangle = 2$ for all~$\alpha \in R$. For more information on f\/inite ref\/lection groups we refer the reader to~\cite{Humph}.

A multiplicity function $k$ on the root system $R$ is a $G$-invariant function $k: R \rightarrow \mC$, i.e.\ $k(\alpha) = k(h \alpha)$ for all $h \in G$. We will denote $k(\alpha)$ by $k_{\alpha}$.

Fixing a positive subsystem $R_{+}$ of the root system $R$ and a multiplicity function $k$, we introduce the Dunkl operators $T_{i}$ associated to $R_{+}$ and $k$ by (see~\cite{MR951883, MR1827871})
\[
T_{i} f(x)= \partial_{x_{i}} f(x) + \sum_{\alpha \in R_{+}} k_{\alpha} \alpha_{i} \frac{f(x) - f(r_{\alpha}(x))}{\langle \alpha, x\rangle}, \qquad f \in C^{1}(\mR^{m}).
\]
An important property of the Dunkl operators is that they commute, i.e.\ $T_{i} T_{j} = T_{j} T_{i}$.

The Dunkl Laplacian is given by $\Delta_{k} = \sum\limits_{i=1}^{m} T_i^2$, or more explicitly by
\[
\Delta_{k} f(x) = \Delta f(x) + 2 \sum_{\alpha \in R_{+}} k_{\alpha} \left( \frac{\langle \nabla f(x), \alpha \rangle}{\langle \alpha, x \rangle}  - \frac{f(x) - f(r_{\alpha}(x))}{\langle \alpha, x \rangle^{2}} \right)
\]
with $\Delta$ the classical Laplacian and $\nabla$ the gradient operator.


If we let $\Delta_{k}$ act on $|x|^2$ we f\/ind $\Delta_{k} |x|^2 = 2m + 4 \gamma = 2 \mu$, where $\gamma = \sum\limits_{\alpha \in R_+} k_{\alpha}$. We call $\mu$ the Dunkl dimension, because most special functions related to $\Delta_{k}$ behave as if one would be working with the classical Laplace operator in a space with dimension $\mu$. We also denote by $\cH_k$ the space of Dunkl-harmonics of degree $k$, i.e.\ $\cH_k = \cP_k \cap \ker{\Delta_{k}}$. The space of Dunkl-harmonics of degree $k$ has the same dimension as the classical space of spherical harmonics of degree $k$ and a basis can e.g.\ be constructed using Maxwell's representation (see \cite{Xu}).


The operators
\[ E:= \tfrac{1}{2}|x|^2,\qquad F:= -\tfrac{1}{2}\Delta_{k} \qquad\text{and}\qquad
H:= \mE +\mu/2
\]
on $\mathcal P$ again satisfy the def\/ining relations of the Lie algebra $\mathfrak{sl}_2$ (see e.g.~\cite{He}).
They are given by
\begin{gather}
\label{sl2reldunkl}
\big[H,E\big] = 2E,\qquad \big[H,F\big] = -2F,\qquad \big[E,F\big] = H.
\end{gather}
As a consequence of these commutation relations, we have the following lemma, which can be proven using induction.

\begin{lemma}\label{calculusDelta}
Let $s \in \mN$ and $R_k \in \cP_k$, then
\[
\Delta_{k} (|x|^{2s} R_{k})= 2s(2k+\mu+2s-2) |x|^{2s-2} R_k + |x|^{2s} \Delta_{k} R_k.
\]
\end{lemma}

\begin{proof}
See \cite{MR1827871}, Lemma 5.1.9.
\end{proof}

This lemma allows us to prove the so-called Fischer decomposition for the Dunkl Laplacian (see e.g.~\cite{MR2207700, MR1827871}).

\begin{theorem}\label{FischerDecomp}
If $\mu \not \in -2\mN$, the space $\cP_k$ decomposes as
\[
\cP_k = \bigoplus_{i=0}^{\lfloor \frac{k}{2}\rfloor} |x|^{2i} \cH_{k-2i}.
\]
\end{theorem}

If we introduce the Dunkl version of the Laplace--Beltrami operator by $\Delta_{\rm LB} = |x|^{2} \Delta_{k} - \mE(\mu - 2 + \mE)$, we can construct projection operators on the dif\/ferent summands in the Fischer decomposition. Indeed, as $\Delta_{\rm LB}$ commutes with $|x|^{2}$ and as
\[
\Delta_{\rm LB}|x|^{2i} \cH_{k-2i} = -(k-2i)(\mu -2 + k - 2i)|x|^{2i} \cH_{k-2i},
\]
where the eigenvalue $(k-2i)(\mu -2 + k - 2i)$ is dif\/ferent for all values of $i$, we immediately have that the operator
\begin{gather}
\mP_i^k = \prod_{l=0, \;  l \neq i}^{\left\lfloor \frac{k}{2} \right\rfloor} \dfrac{\Delta_{\rm LB} + (k-2l)(\mu-2+k-2l)}{2(i-l)(2k-2i-2l+\mu-2)}
\label{projoperators}
\end{gather}
satisf\/ies
\[
\mP_i^k (|x|^{2j} \cH_{k-2j}) = \delta_{ij} |x|^{2i} \cH_{k-2i}.
\]
This provides an easier way to construct projection operators than the method presented in \cite{MR1827871} and \cite{MR2207700}.

Now let $\{\phi_\nu\,, \nu\in \mZ_+^m\}$ be a basis of $\mathcal P$ such that $\phi_\nu\in
\mathcal P_{|\nu|}$. The Hermite polynomials related to $G$ are def\/ined as follows by R\"osler (see \cite{MR1620515}).

\begin{definition}\label{DefHermRos}
The generalized Hermite polynomials $\{H_\nu\,, \>\nu\in \mZ_+^m\}$
associated with the basis~$\{\phi_\nu\}$ on $\mR^m$ are given by
\begin{gather*}
H_\nu(x):= 2^{|\nu|}e^{-\Delta_{k}/4}\phi_\nu(x) =
2^{|\nu|}\sum_{n=0}^{\lfloor|\nu|/2\rfloor} \frac{(-1)^n}{4^n n!}\,
\Delta_{k}^n \phi_\nu(x).
\end{gather*}
Moreover, the generalized Hermite functions on $\mR^m$ are def\/ined by
\begin{gather*}
h_\nu(x):= e^{-|x|^2/2}H_\nu(x), \qquad \nu\in \mZ_+^m.
\end{gather*}
\end{definition}

In \cite{MR1620515}, for some specif\/ic ref\/lection groups $G$, these Hermite polynomials are expressed in terms of known special functions (such as the Jack polynomials). However, their general structure is not very clear from Def\/inition~\ref{DefHermRos}.

In the sequel, we will also need the following theorem (see \cite{MR1620515}).
\begin{theorem}\quad{}\label{DiffEqRosler}
\begin{enumerate}\itemsep=0pt
\item[\rm{(1)}] For $n\in \mZ_+$ set $ V_n:= \{e^{-\Delta_{k}/4} p :
p\in \mathcal P_n\}.$ Then $ \mathcal P = \bigoplus_{n\in \mZ_+} V_n$,
and $V_n$ is the eigenspace of the operator $ \Delta_{k} -2\mE $ on
$\mathcal P$ corresponding to the eigenvalue $-2n$.
\item[\rm{(2)}] For $q\in V_n$, the function $ f(x):= e^{-|x|^2/2} q(x)$
satisfies
\[ \bigl(\Delta_{k} -|x|^2\bigr) f = -(2n + \mu)f .\]
\end{enumerate}
\end{theorem}

The Dunkl transform $\cD$ (see \cite{Du4, deJ}) associated with $G$ and $k\geq 0$ is def\/ined by
\begin{gather*}
 \mathcal D: \ \ L^1(\mR^m, w_k(x)dx)\to C(\mR^m);\qquad
 \mathcal D f(\xi):= \int_{\mR^m} f(x)\,K(-i\xi,x)\,w_k(x)dx \quad
(\xi\in \mR^m),
\end{gather*}
where $K(-i\xi,x)$ is the so-called Dunkl kernel or generalized exponential and with $w_k(x) = \prod\limits_{\alpha \in R_{+}} |\langle \alpha, x\rangle |^{2 k_{\alpha}}$ the weight function corresponding to $G$.

The following proposition is proven in \cite{MR1620515} and gives the action of the Dunkl transform on the generalized Hermite polynomials.
\begin{proposition}\label{FourierHermite}
The generalized Hermite functions $\{h_\nu,\, \nu\in
\mZ_+^m\}$ are a basis of eigenfunctions of the Dunkl transform
$\mathcal D$ on $L^2(\mR^m, w_k(x)dx)$, satisfying
\[
\mathcal D(h_\nu)=  2^{\mu/2} c_k^{-1} (-i)^{|\nu|} h_\nu
\]
with $2^{\mu/2} c_k^{-1} = \cD(e^{-|x|^2/2})(0)$.
\end{proposition}

Note that it is possible to introduce a Dunkl version of the Dirac operator. Indeed, as the Dunkl operators $T_{i}$ are commutative, we can factorize the Dunkl Laplacian in the same way as the usual Laplace operator. Def\/ining the Dunkl Dirac operator~$D_{k}$ by
\[
D_{k} = \sum_{i=1}^{m} e_{i} T_{i}
\]
with the $e_{i}$ generators of the orthogonal Clif\/ford algebra as in Section~\ref{hermCA}, we obtain that $D_{k}^{2} = - \Delta_{k}$. For some basic results on the Dunkl Dirac operator, such as a Cauchy formula, we refer the reader to~\cite{MR2230262}.

For the sequel we need the anti-commutator of $D_{k}$ and $\ux$. We obtain
\begin{gather*}
\left\{ D_{k}, \ux \right\} = \sum_{i,j} \{e_{i} T_{i}, e_{j} x_{j}\}
= \sum_{i,j} (e_{i} e_{j} T_{i} x_{j} + e_{j} e_{i} x_{j} T_{i})\\
\phantom{\left\{ D_{k}, \ux \right\}}{} = -\sum_{i} (T_{i} x_{i} + x_{i} T_{i}) + \sum_{i \neq j} e_{i} e_{j} (T_{i} x_{j} - x_{j} T_{i}).
\end{gather*}
The f\/irst term $\sum_{i} (T_{i} x_{i} + x_{i} T_{i})$ equals $2\mE + \mu$ (see e.g.~\cite{Said}, formulae~(3.5) and~(1.2)).

For the second term, we calculate the action of $T_{i} x_{j} - x_{j} T_{i}$ on a function $f$, yielding
\begin{gather*}
(T_{i} x_{j}  - x_{j} T_{i}) f = \sum_{\alpha \in R_{+}} k_{\alpha} \alpha_{i}\frac{x_{j} f(x) - (r_{\alpha}(x))_{j} f(r_{\alpha}(x)) }{\langle \alpha, x\rangle} - x_{j}\sum_{\alpha \in R_{+}} k_{\alpha} \alpha_{i}\frac{ f(x) - f(r_{\alpha}(x)) }{\langle \alpha, x\rangle}\\
\phantom{(T_{i} x_{j}  - x_{j} T_{i}) f}{} = \sum_{\alpha \in R_{+}} k_{\alpha} \alpha_{i}\frac{ \langle \alpha,x\rangle \alpha_{j} f(r_{\alpha}(x)) }{\langle \alpha, x\rangle}
= \sum_{\alpha \in R_{+}} k_{\alpha} \alpha_{i} \alpha_{j} f(r_{\alpha}(x)),
\end{gather*}
an expression symmetric in $i$ and $j$. As $e_{i} e_{j}$ is anti-symmetric in $i$ and $j$, the sum $\sum\limits_{i \neq j} e_{i} e_{j} (T_{i} x_{j} - x_{j} T_{i})$ vanishes.
We conclude that
\begin{gather}
\left\{ D_{k}, \ux \right\} = -(2\mE + \mu).
\label{anticommDiracX}
\end{gather}

\section[Clifford-Hermite polynomials related to the Dunkl Laplacian]{Clif\/ford--Hermite polynomials related to the Dunkl Laplacian}
\label{CHdunkl}

In this section we construct the Clif\/ford--Hermite polynomials related to the Dunkl Laplacian. By analogy with Section~\ref{hermCA}, we f\/irst introduce the operator
\[
D_{+} =  - D_{k} + 2 \ux.
\]
Calculating the square of this operator yields
\begin{gather*}
D_{+}^{2} = (-D_{k}+2\ux)^{2}
= -\Delta_{k} + 4 \ux^{2} - 2\left\{ D_{k}, \ux \right\}\\
\phantom{D_{+}^{2}}{} =-\Delta_{k} - 4 |x|^2 + 2(2\mE + \mu),
\end{gather*}
where we have used formula (\ref{anticommDiracX}).

We then def\/ine the Clif\/ford--Hermite polynomials as follows.

\begin{definition}
Let $H_k \in \cH_k$ be a Dunkl-harmonic of degree $k$ and $t$ a positive integer. Then
\[
C\!H_{2t}^{\mu}(H_k) = (D_+)^{2t} H_k
\]
is a Clif\/ford--Hermite polynomial of degree $2t$ associated with $H_{k}$.
\end{definition}

Using Lemma~\ref{calculusDelta} we immediately see that the precise form of the polynomials $C\!H_{2t}^{\mu}(H_k)$ depends only on the degree of the spherical harmonic $H_k$, so we can write
\[
C\!H_{2t}^{\mu}(H_k)=\phi_{2t,k}^{\mu}(|x|^{2})H_k.
\]
The f\/irst few Clif\/ford--Hermite polynomials have the following explicit form:
\begin{gather*}
C\!H_{0}^{\mu}(H_k) =H_k,\\
C\!H_{2}^{\mu}(H_k) = \left( -4 |x|^2 + 2(2k+\mu)\right) H_k,\\
C\!H_{4}^{\mu}(H_k) =\left( 16 |x|^4 - 16 (2k+\mu+2)|x|^2+4(2k+\mu+2)(2k+\mu)\right) H_k.
\end{gather*}

Using the def\/inition, we immediately obtain the following recursion relation
\begin{gather}
C\!H_{2t}^{\mu}(H_k) = (D_+)^2 C\!H_{2t-2}^{\mu}(H_k).
\label{recursionCH}
\end{gather}

There also exists a Rodrigues formula for the Clif\/ford--Hermite polynomials.

\begin{theorem}[Rodrigues formula] The Clifford--Hermite polynomials take the form
\begin{gather*}
C\!H_{2t}^{\mu}(H_k) = \exp(|x|^2/2) (-\Delta_{k} - |x|^2 + 2 \mE + \mu)^{t} \exp(-|x|^2/2) H_k\\
\phantom{C\!H_{2t}^{\mu}(H_k)}{} = \exp(|x|^2) (-\Delta_{k})^t  \exp(-|x|^2) H_k.
\end{gather*}
\end{theorem}

\begin{proof}
This follows immediately from the following operator equalities:
\[
-\exp(|x|^2) \Delta_{k}  \exp(-|x|^2) = (D_+)^2 = \exp(|x|^2/2) (-\Delta_{k} -|x|^2 + 2 \mE + \mu) \exp(-|x|^2/2),
\]
which can e.g.~be found in \cite[pp. 254--255]{Said}, combined with the def\/inition of the Clif\/ford--Hermite polynomials.
\end{proof}

Now we prove that the Clif\/ford--Hermite polynomials satisfy a partial dif\/ferential equation.

\begin{theorem}[Dif\/ferential equation]\label{DiffEqCH}
$C\!H_{2t}^{\mu}(H_k)$ is a solution of the following differential equation:
\[
[\Delta_{k} - 2 \mE ] C\!H_{2t}^{\mu}(H_k) = -2(2t+k) C\!H_{2t}^{\mu}(H_k).
\]
\end{theorem}

\begin{proof}
We can expand the Clif\/ford--Hermite polynomials as follows:
\[
C\!H_{2t}^{\mu}(H_k)= \sum_{i=0}^{t}a_{2i}^{2t} |x|^{2i} H_k.
\]
The recursion relation (\ref{recursionCH}) leads to the following relation between the coef\/f\/icients $a_{2i}^{2t}$:
\begin{gather}
a_{2i}^{2t} = -(2i+2)(2k+\mu+2i)a_{2i+2}^{2t-2} + 2(2k + 4i +\mu) a_{2i}^{2t-2} - 4 a_{2i-2}^{2t-2}.
\label{RecursionFormula}
\end{gather}
We need to prove that
\[
 (2i+2)(2k+\mu+2i) a_{2i+2}^{2t} - 2(2i + k) a_{2i}^{2t}  = -2(2t+k) a_{2i}^{2t}
\]
or
\begin{gather}
- 2(2t-2i) a_{2i}^{2t} = (2i+2)(2k+\mu+2i) a_{2i+2}^{2t}.
\label{DiffeqFormula}
\end{gather}
It is easy to check that the theorem holds for $t=0$. Substituting (\ref{RecursionFormula}) in (\ref{DiffeqFormula}) and using induction on $t$ then completes the proof.
\end{proof}

Using the proof of the previous lemma (formulae (\ref{DiffeqFormula}) and (\ref{RecursionFormula})) it is now possible to obtain the following closed form for the Clif\/ford--Hermite polynomials.
\begin{theorem}\label{CHLagRelation}
The Clifford--Hermite polynomials can be written in terms of the generalized Laguerre polynomials as
\begin{gather*}
\phi_{2t,k}^{\mu}(|x|^{2}) = 2^{2t} t! L_{t}^{\frac{\mu}{2} + k-1}(|x|^2),
\end{gather*}
with
\[
L_{t}^{\alpha}(x) = \sum_{i=0}^t \frac{\Gamma(t +\alpha +1)}{i! (t-i)! \Gamma(i + \alpha +1)} (-x)^i.
\]
\end{theorem}

\begin{proof}
Recall that
\[
\phi_{2t,k}^{\mu}(|x|^{2})= \sum_{i=0}^{t}a_{2i}^{2t} |x|^{2i}.
\]
Using formula (\ref{DiffeqFormula}) we obtain
\begin{gather*}
a_{2i}^{2t} = - \frac{t-i+1}{i (k + \mu/2 +i-1 )} a_{2i-2}^{2t} =\cdots  = (-1)^{i} \frac{t!}{i! (t-i)!} \frac{\Gamma(k + \mu/2)}{\Gamma(k + \mu/2 +i)} a_{0}^{2t}.
\end{gather*}
Using formula (\ref{RecursionFormula}) and again (\ref{DiffeqFormula}) we also have
\begin{gather*}
a_{0}^{2t} = 2 (2k + \mu) \big[a_{0}^{2t-2} - a_{2}^{2t-2}\big]
= 4 (k + \mu/2 + t-1) a_{0}^{2t-2} = \cdots\\
\phantom{a_{0}^{2t}}{}= 2^{2t} \frac{\Gamma(k + \mu/2 + t)}{\Gamma(k + \mu/2)}a_{0}^{0}
= 2^{2t} \frac{\Gamma(k + \mu/2 + t)}{\Gamma(k + \mu/2)}.
\end{gather*}
Combining both results and comparing with the def\/inition of the generalized Laguerre polynomials yields the result of the theorem.
\end{proof}

Using Theorem \ref{CHLagRelation} and multiplying the Clif\/ford--Hermite polynomials with the Gaussian we obtain the set of functions
\[
L_{t}^{\frac{\mu}{2} + k-1}(|x|^2) H_{k} e^{-|x|^{2}/2}, \qquad H_{k} \in \cH_{k}.
\]
These functions have previously also been introduced by Dunkl in~\cite{Du4} as a basis of the weighted $L^{2}$ space $L^{2}(\mR^{m}, w_{k}(x) dx)$. Recently they have been studied more thorougly in~\cite{Said} to prove that the $\mathfrak{sl}_{2}$ relations (\ref{sl2reldunkl}) exponentiate to a unique unitary representation of the universal covering group of $SL(2, \mR)$.

Now we are able to state the connection between the Clif\/ford--Hermite polynomials and the generalized Hermite polynomials of R\"osler. First we def\/ine the subspace $W_n \subset \cP$ by
\[
W_n = \bigoplus_{k+ 2t = n} \phi_{2t,k}^{\mu}(|x|^{2})\cH_k.
\]
Note that, because of Theorem~\ref{FischerDecomp}, we have that $\oplus_{n=0}^{\infty}W_n = \cP$. Also, each summand $\phi_{2t,k}^{\mu}(|x|^{2})\cH_k$ is clearly invariant under the action of $G$.

We then have the following theorem.
\begin{theorem}
The spaces $W_n$ and $V_n$, defined in Theorem~{\rm \ref{DiffEqRosler}}, coincide.
\label{equalityHerm}
\end{theorem}

\begin{proof}
Both $W_n$ and $V_n$ are maximal eigenspaces in $\cP$ for the operator $\Delta_{k} - 2\mE$, corresponding with the eigenvalue $-2n$ (see Theorems \ref{DiffEqRosler} and \ref{DiffEqCH}). Hence, we have that $W_n = V_n$.
\end{proof}

This theorem allows us to make the link between the Clif\/ford--Hermite and the generalized Hermite polynomials even more explicit. If we choose a basis $\{\psi_{j}\}$ of $\cP_{n}$ in such a way that each basis element is of the form
\[
|x|^{2i} H_{n-2i}, \qquad H_{n-2i} \in \cH_{n-2i},
\]
which is always possible due to Theorem~\ref{FischerDecomp}, then the generalized Hermite polynomial
\[
2^{n}e^{-\Delta_{k}/4} \left( |x|^{2i} H_{n-2i} \right)
\]
will be proportional to the Clif\/ford--Hermite polynomial
\[
\phi_{2i,k}^{\mu}(|x|^{2}) H_{n-2i}.
\]

Theorem \ref{equalityHerm} allows us also to transfer results from \cite{MR1620515} to the Clif\/ford--Hermite polynomials and vice versa. As an example, using Proposition~\ref{FourierHermite} we obtain that every function of the form
\[
\phi_{2t,k}^{\mu}(|x|^{2}) H_{k} e^{-|x|^2/2}, \qquad H_{k} \in \cH_{k}
\]
is an eigenfunction of the Dunkl transform with corresponding eigenvalue
\[
2^{\mu/2} c_k^{-1}\,(-i)^{2t + k},
\]
a result which can also be found in \cite{Du4, Said}.

\begin{remark}
Note that each summand $\phi_{2t,k}^{\mu}(|x|^{2})\cH_k$ in the decomposition of $W_{n}$ is an eigenspace of the Laplace--Beltrami operator $\Delta_{\rm LB}$ with eigenvalue $-k(\mu -2 + k)$. Hence it is possible to construct projection operators on each summand in the decomposition of $W_{n}$ in a similar way as in formula (\ref{projoperators}).
\end{remark}

\begin{remark}
Note that it is also possible to study generalized Gegenbauer polynomials with respect to the Dunkl Laplacian. This can be done in a similar way as in \cite{DBS3} for the case of $O(m) \times Sp(2n)$.
\end{remark}

\subsection*{Acknowledgements}

The author is supported by a Ph.D. Fellowship of the the Research Foundation - Flanders (FWO).

\pdfbookmark[1]{References}{ref}
\LastPageEnding

\end{document}